\documentclass[11pt,reqno]{amsart}
\usepackage{amssymb,amscd,amsbsy}
\usepackage{amssymb,amscd,amsbsy,mathrsfs}
\newcommand{\lb}{\linebreak}

\renewcommand{\a}{\alpha}
\renewcommand{\b}{\beta}

\renewcommand{\d}{\delta}

\newcommand{\vk}{\varkappa}
\newcommand{\z}{\zeta}

\newcommand{\s}{\sigma}
\renewcommand{\t}{\tau}

\newcommand{\f}{\varphi}
\renewcommand{\o}{\omega}

\newcommand{\D}{\Delta}
\renewcommand{\L}{\Lambda}

\newcommand{\B}{{\mathscr B}}

\newcommand{\E}{{\mathscr E}}
\newcommand{\cd}{{\mathscr D}}
\newcommand{\F}{{\mathscr F}}

\newcommand{\h}{{\mathscr H}}

\newcommand{\X}{{\mathscr X}}
\newcommand{\Y}{{\mathscr Y}}

\newcommand{\C}{{\Bbb C}}

\newcommand{\R}{{\Bbb R}}
\newcommand{\Z}{{\Bbb Z}}

\newcommand{\0}{{\boldsymbol{0}}}

\newcommand{\bs}{\boldsymbol}

\newcommand{\bS}{{\boldsymbol S}}

\newcommand{\rf}[1]{(\ref{#1})}

\newcommand{\df}{\stackrel{\mathrm{def}}{=}}

\newcommand{\re}{\operatorname{Re}}

\newcommand{\supp}{\operatorname{supp}}

\newcommand{\trace}{\operatorname{trace}}
\newcommand{\rank}{\operatorname{rank}}
\newcommand{\const}{\operatorname{const}}

\newcommand{\eeq}{\end{equation}}
\newcommand{\beq}{\begin{equation}}
\newcommand{\bay}{\begin{eqnarray}}
\newcommand{\ba}{\begin{align*}}
\newcommand{\ea}{\end{align*}}
\newcommand{\ey}{\end{eqnarray}}
\newcommand{\bey}{\begin{eqnarray*}}
\newcommand{\eey}{\end{eqnarray*}}

\newcommand{\be}{\infty}

\newcommand{\bl}{\blacksquare}

\newcommand{\Pf}{{\bf Proof. }}
\newcommand{\im}{\operatorname{Im}}
\renewcommand{\re}{\operatorname{Re}}

\newtheorem{thm}{\hspace{\parindent}Theorem}[section]

\newtheorem{cor}[thm]{\hspace{\parindent}Corollary}
\newtheorem{lem}[thm]{\hspace{\parindent}Lemma}

\theoremstyle{remark}

\newtheorem*{rem*}{Remark}

\newcommand\Li{{\rm Lip}}
\newcommand\fM{\frak M}

\newcommand\dg{\frak D}

\newcommand{\fI}{{\frak I}}

\newcommand\mC{\mathcal{C}}
\newcommand\mQ{\mathcal{Q}}
\newcommand\mR{\mathcal{R}}

\newcommand{\qm}{\quad\mbox{and}\quad}

\begin{document}

\newcommand{\vse}{\vspace{.2in}}
\numberwithin{equation}{section}

\title{Functions of perturbed $\bs{n}$-tuples of commuting self-adjoint operators}
\author{F.L. Nazarov and V.V. Peller}

\begin{abstract}
Let $(A_1,\cdots,A_n)$ and $(B_1,\cdots,B_n)$ be $n$-tuples of commuting self-adjoint operators on Hilbert space. For functions $f$ on $\R^n$ satisfying certain conditions, we obtain sharp estimates of the operator norms (or norms in operator ideals) of
$f(A_1,\cdots,A_n)-f(B_1,\cdots,B_n)$ in terms of the corresponding norms of $A_j-B_j$, $1\le j\le n$.

We obtain analogs of earlier results on estimates for functions of perturbed self-adjoint and normal operators. It turns out that for $n\ge3$, the methods that were used for self-adjoint and normal operators do not work. We propose a new method that works for arbitrary $n$.

We also get sharp estimates for quasicommutators
$f(A_1,\cdots,A_n)R-Rf(B_1,\cdots,B_n)$ in terms of norms of
$A_jR-RB_j$, $1\le j\le n$, for a bounded linear operator $R$.
\end{abstract}

\maketitle

\

\begin{center}
{\Large Contents}
\end{center}

\

\begin{enumerate}
\item[1.] Introduction \quad\dotfill \pageref{In}
\item[2.] Preliminaries \quad\dotfill \pageref{prel}
\item[3.] The techniques used in the case of normal operators do not work when \lb $n\ge3$
\quad\dotfill \pageref{n=3}
\item[4.] An integral formula  \quad\dotfill \pageref{intf}
\item[5.] Operator Lipschitzness \quad\dotfill \pageref{Olss}
\item[6.] H\"older classes and arbitrary moduli of continuity \quad\dotfill \pageref{HcMc}
\item[7.] Estimates in ideal norms \quad\dotfill \pageref{Idealc}
\item[8.] Commutator and quasicommutator estimates \quad\dotfill \pageref{ComQ}
\item[] References \quad\dotfill \pageref{bibl}
\end{enumerate}

\

\setcounter{section}{0}
\section{\bf Introduction}
\setcounter{equation}{0}
\label{In}

\medskip

In this paper we extend earlier results of \cite{Pe2}, \cite{Pe3}, \cite{AP2}, \cite{AP3}, \cite{AP4}, \cite{APPS2} on functions of perturbed self-adjoint and normal operators to the case of functions of perturbed $n$-tuples of commuting self-adjoint operators.

It is well known that a Lipschitz function $f$ on the real line ${\Bbb R}$ (i.e., a function $f$ satisfying the inequality $|f(x)-f(y)|\le\const|x-y|$, \:\mbox{$x,y\in\R$}) does not have to satisfy the inequality
\bay
\label{OpL}
\|f(A)-f(B)\|\le\const\|A-B\|
\ey
for all bounded self-adjoint operators $A$ and $B$ on Hilbert space. 
In other words, it does not have to be {\it operator Lipschitz}.
Note that if $f$ is
operator Lipschitz, then inequality \rf{OpL} also holds for unbounded self-adjoint operators
$A$ and $B$ as far as $A-B$ is bounded, see \cite{AP2} and \cite{AP4}. 

The first example of a Lipschitz
function that is not operator Lipschitz was obtained by Farforovskaya \cite{F}. We refer the
reader to \cite{Pe2} (see also \cite{Pe3}) where necessary conditions and sufficient conditions for a function on $\R$ to be operator Lipschitz are obtained. In particular, it was shown in \cite{Pe2} and \cite{Pe3} that functions in the Besov space $B_{\be,1}^1(\R)$ are operator Lipschitz (see Subsection 2.1 below for a brief introduction to Besov spaces). On the other hand, operator Lipschitz functions must belong locally to $B_{1,1}^1(\R)$, see \cite{Pe2} and \cite{Pe3}.

It turns out, however, that the situation is entirely different if we proceed from Lipschitz functions to H\"older functions. It was shown in \cite{AP2} (see also \cite{AP1}) that if $f$ belongs to the H\"older class $\L_\a(\R)$, $0<\a<1$, i.e., $|f(x)-f(y)|\le\const|x-y|^\a$, $x,\,y\in\R$, then $f$ is necessarily {\it operator H\"older of order} $\a$, i.e.,
$$
\|f(A)-f(B)\|\le\const\|A-B\|^\a
$$
for all self-adjoint operators $A$ and $B$ on Hilbert space with bounded $A-B$. 
Note that in \cite{AP1} and \cite{AP2} sharp results were also obtained for functions in the space $\L_\o$ for an arbitrary modulus of continuity $\o$.
Similar (slightly weaker) results were obtained independently in \cite{FN}.  

In this paper we are also going to consider the case of perturbations by operators of Schatten--von Neumann class $\bS_p$. It was proved in \cite{AP3} (see also \cite{AP1}) that for $f\in\L_\a(\R)$, $0<\a<1$, $p>1$, and for self-adjoint operators $A$ and $B$ with $A-B\in\bS_p$, the operator $f(A)-f(B)$ must be in $\bS_{p/\a}$ and the following inequality holds:
$$
\|f(A)-f(B)\|_{\bS_{p/\a}}\le\const\|A-B\|_{\bS_p}^\a.
$$
Let us also mention that in \cite{AP3} more general results for operator ideal (quasi)norms were obtained as well.

It turns out that the situation for functions of normal operators or, which is the same, for functions of two commuting self-adjoint operators is more complicated and requires different techniques. Nevertheless in \cite{APPS2} (see also \cite{APPS}) analogs of the above mentioned results are obtained for normal operators and functions on the plane. 

We consider the more general case of functions of $n$-tuples of commuting self-adjoint operators and find sharp estimates of $f(A_1,\cdots,A_n)-f(B_1,\cdots,B_n)$ for functions $f$ on $\R^n$.
It turns out that the techniques used in \cite{APPS2} to obtain results for normal operators do not work in the case $n\ge3$. 

We propose a different method, which allows us to obtain analogs of the above results for functions of $n$-tuples of commuting self-adjoint operators for arbitrary $n$. Our results considerably improve the results of \cite{F2}.

Our results are based on the {\it crucial estimate} obtained in Theorem 5.1 (see \S\,\ref{Olss}):
$$
\big\|f(A_1,\cdots,A_n)-f(B_1,\cdots,B_n)\big\|
\le c_n\s\max_{1\le j\le n}\|A_j-B_j\|
$$
for some positive number $c_n$, whenever $f$ is a bounded function on $\R^n$ such that its Fourier transform $\F f$ is supported in $\{\z\in\R^n:~|\xi|\le\s\}$, and $(A_1,\cdots,A_n)$ and $(B_1,\cdots,B_n)$ are $n$-tuples of commuting self-adjoint operators.

In \S\,3 we explain why the methods used in \cite{APPS2} for normal operators do not work for $n$-tuples of commuting self-adjoint operators if $n\ge3$.

We establish in \S\,4 a formula that, for $n$-tuples $(A_1,\cdots,A_n)$ and $(B_1,\cdots,B_n)$ of commuting self-adjoint operators and suitable functions $f$ on $\R^n$,
represents the difference $f(A_1,\cdots,A_n)-f(A_1,\cdots,A_n)$ in terms of double operator integrals.

We obtain in \S\,5 sharp sufficient conditions on functions on $\R^n$ to be operator Lipschitz and obtain Lipschitz type estimates in operator ideals.

In \S\,6 we prove that H\"older functions of order $\a$, $0<\a<1$, on $\R^n$ are   
operator H\"older of order $\a$. We also consider the case of arbitrari moduli of continuity. 

Section 7 is devoted to ideal (quasi)norm estimates of $f(A_1,\cdots,A_n)-f(B_1,\cdots,B_n)$ for H\"older functions $f$ of order $\a$, $0<\a<1$.

Finally, in \S\,8 we obtain quasicommutator estimates, i.e., estimates of
the operators $f(A_1,\cdots,A_n)R-Rf(B_1,\cdots,B_n)$ in terms of norms of the quasicommutators $A_jR-RA_j$, $1\le j\le n$, where $R$ is a bounded linear operator.

We collect in \S\,2 necessary information on function spaces, operator ideals and double operator integrals.

Note that main results of this paper were announced in the note \cite{NP}.

We are grateful to A.B. Aleksandrov for helpful discussions.

\

\section{\bf Preliminaries}
\setcounter{equation}{0}
\label{prel}

\

In this section we collect necessary information on function spaces, operator ideals, and double operator integrals. 

\medskip

{\bf 2.1. Littelewood--Paley type expansions and function spaces.} The technique of Littlewood--Paley type expansions of functions or distributions on Euclidean spaces 
is a very important tool in Harmonic analysis. 

Let $w$ be an infinitely differentiable function on $\R$ such
that
\bay
\label{w}
w\ge0,\quad\supp w\subset\left[\frac12,2\right],\quad\mbox{and} \quad w(s)=1-w\left(\frac s2\right)\quad\mbox{for}\quad s\in[1,2].
\ey

We define the functions $W_l$, $l\in\Z$, on $\R^n$ by 
$$
\big(\F W_l\big)(x)=w\left(\frac{|x|}{2^l}\right),\quad l\in\Z, \quad x=(x_1,\cdots,x_n),
\quad|x|\df\left(\sum_{j=1}^nx_j^2\right)^{1/2},
$$
where $\F$ is the {\it Fourier transform} defined on $L^1\big(\R^n\big)$ by
$$
\big(\F f\big)(t)=\!\int\limits_{\R^n} f(x)e^{-{\rm i}(x,t)}\,dx,\!\quad 
x=(x_1,\cdots,x_n),
\quad t=(t_1,\cdots,t_n), \!\quad(x,t)\df \sum_{j=1}^nx_jt_j.
$$
Clearly,
$$
\sum_{l\in\Z}(\F W_l)(t)=1,\quad t\in\R^n\setminus\{0\}.
$$

With each tempered distribution $f\in{\mathscr S}^\prime\big(\R^n\big)$, we
associate the sequence $\{f_l\}_{l\in\Z}$,
\bay
\label{fn}
f_l\df f*W_l.
\ey
The formal series
$$
\sum_{l\in\Z}f_l
$$
is a Littlewood--Paley type expansion of $f$. This series does not necessarily converge to $f$. However, for certain function spaces we deal with in this paper we have equality
$$
f(x)-f(y)=\sum_{l\in\Z}\big(f_l(x)-f_l(y)\big),\quad x,~y\in\R^n,
$$
and the series on the right converges uniformly. This allows us to use the formula
$$
f(A_1,\cdots,A_n)-f(B_1,\cdots,B_n)=
\sum_{l\in\Z}\big(f_l(A_1,\cdots,A_n)-f_l(B_1,\cdots,B_n)\big)
$$
for $n$-tuples of commuting self-adjoint operators $(A_1,\cdots,A_n)$ and 
$(B_1,\cdots,B_n)$.

In this paper an important role is played by the spaces $\L_\o(\R^n)$, where $\o$ is a {\it modulus of continuity}, i.e., $\o$ is a nondecreasing continuous function on $[0,\be)$
such that $\o(0)=0$, $\o(x)>0$ for $x>0$, and
$$
\o(x+y)\le\o(x)+\o(y),\quad x,~y\in[0,\be).
$$
The space $\L_\o\big(\R^n\big)$ consists of functions $f$ on $\R^n$ such that
$$
\|f\|_{\L_\o}\df\sup_{x\ne y}\frac{|f(x)-f(y)|}{\o(|x-y|)}<\be.
$$
Functions in $\L_\o\big(\R^n\big)$ satisfy the following Bernstein type estimates:
$$
f\in \L_\o\big(\R^n\big)\quad\Longrightarrow\quad
\left\|f-\sum_{l=m}^\be f*W_l\right\|_{L^\be}\le\const\o\big(2^{-m}\big)\|f\|_{\L_\o},\quad m\in\Z,
$$
and so
$$
f\in \L_\o\big(\R^n\big)\quad\Longrightarrow\quad
\|f*W_l\|_{L^\be}\le\const\o\big(2^{-l}\big)\|f\|_{\L_\o},\quad l\in\Z,
$$
(see \cite{AP2} and \cite{APPS2} for details).

In the case $\o(t)=t^\a$, $0<\a<1$, we use the notation 
$\L_\a\big(\R^n\big)\df\L_\o\big(\R^n\big)$. The space $\L_\a\big(\R^n\big)$ is called the {\it space of H\"older functions of order} $\a$. It admits the following characterization:
$$
f\in \L_\a\big(\R^n\big)\quad\Longleftrightarrow\quad\|f*W_l\|_{L^\be}\le\const2^{-l\a},\quad l\in\Z.
$$

The H\"older classes $\L_\a\big(\R^n\big)$ form a special case of Besov spaces that play an important role in problems of perturbation theory.

Initially we define the (homogeneous) Besov class $\dot B^s_{pq}\big(\R^n\big)$,
$s>0$, $1\le p,\,q\le\be$, as the space of all
$f\in{\mathscr S}^\prime(\R^n)$
such that
\bay
\label{Wn}
\{2^{ls}\|f_l\|_{L^p}\}_{l\in\Z}\in\ell^q(\Z)
\ey
and put
$$
\|f\|_{B^s_{pq}}\df\big\|\{2^{ls}\|f_l\|_{L^p}\}_{l\in\Z}\big\|_{\ell^q(\Z)}.
$$
According to this definition, the space $\dot B^s_{pq}(\R^n)$ contains all polynomials
and all polynomials $f$ satisfy the equality $\|f\|_{B^s_{pq}}=0$. Moreover, the distribution $f$ is determined by the sequence $\{f_l\}_{n\in\Z}$
uniquely up to a polynomial. It is easy to see that the series 
$\sum_{l\ge0}f_l$ converges in ${\mathscr S}^\prime(\R^n)$.
However, the series $\sum_{l<0}f_l$ can diverge in general. It can easily be proved that the series
\bay
\label{ryad}
\sum_{l<0}\frac{\partial^r f_l}{\partial x_1^{r_1}\cdots\partial x_n^{r_n}},\qquad \mbox{where}\quad r_j\ge0,\quad\mbox{for}\quad
1\le j\le n,\quad\sum_{j=1}^nr_j=r,
\ey
converges uniformly on $\R^n$ for every nonnegative integer
$r>s-n/p$. Note that in the case $q=1$ the series \rf{ryad}
converges uniformly, whenever $r\ge s-n/p$.

Now we can define the modified (homogeneous) Besov class $B^s_{pq}\big(\R^n\big)$. We say that a distribution $f$
belongs to $B^s_{pq}(\R^n)$ if \rf{Wn} holds and
$$
\frac{\partial^r f}{\partial x_1^{r_1}\cdots\partial x_n^{r_n}}
=\sum_{l\in\Z}\frac{\partial^r f_l}{\partial x_1^{r_1}\cdots\partial x_n^{r_n}},\quad
\mbox{whenever}\quad 
r_j\ge0,\quad\mbox{for}\quad
1\le j\le n,\quad\sum_{j=1}^nr_j=r.
$$
in the space ${\mathscr S}^\prime\big(\R^n\big)$, where $r$ is
the minimal nonnegative integer such that $r>s-n/p$ ($r\ge s-n/p$ if $q=1$). Now the function $f$ is determined uniquely by the sequence $\{f_l\}_{l\in\Z}$ up
to a polynomial of degree less than $r$, and a polynomial $g$ belongs to $B^s_{pq}\big(\R^n\big)$
if and only if $\deg g<r$.

%To define a regularized de la Vall\'ee Poussin type kernel $V_n$, we define the $C^\be$ function $v$ on $\R$ by
%\bey
%%\label{VP}
%v(x)=1\quad\mbox{for}\quad x\in[-1,1]\quad\mbox{and}\quad v(x)=w(|x|)\quad\mbox{if}\quad |x|\ge1,
%\eey
%where $w$ is the function defined by \rf{w}. Now
%we can define the de la Vall\'ee Poussin type functions $V_n$ by
%$$
%\F V_n(x)=v\left(\frac{|x|}{2^n}\right),\quad n\in\Z,\quad x=(x_1,x_2).
%$$
%We put $V\df V_0$. Clearly, $V_n(x)=2^{2n}V(2^n x)$.

Besov classes admit many other descriptions.
We give here the definition in terms of finite differences.
For $h\in\R^n$, we define the difference operator $\D_h$,
$$
(\D_hf)(x)=f(x+h)-f(x),\quad x\in\R^n.
$$
It is easy to see that $B_{pq}^s\big(\R^n\big)\subset L^1_{\rm loc}\big(\R^n\big)$ for every $s>0$
and $B_{pq}^s\big(\R^n\big)\subset C\big(\R^n\big)$ for every $s>n/p$. Let $s>0$ and let $m$ be the integer such that $m-1\le s<m$.
The Besov space $B_{pq}^s\big(\R^n\big)$ can be defined as the set of
functions $f\in L^1_{\rm loc}\big(\R^n\big)$ such that
$$
\int_{\R^n}|h|^{-n-sq}\|\D^m_h f\|_{L^p}^q\,dh<\be\quad\mbox{for}\quad q<\be
$$
and
\bay
\label{pbe}
\sup_{h\not=0}\frac{\|\D^m_h f\|_{L^p}}{|h|^s}<\be\quad\mbox{for}\quad q=\be.
\ey
However, with this definition the Besov space can contain polynomials of higher degree than in the case of the first definition given above.

We use the notation $B_p^s\big(\R^n\big)$ for $B_{pp}^s\big(\R^n\big)$.

Clearly, $\L_\a\big(\R^n\big)=B^\a_\be\big(\R^n\big)$. Moreover, this equality can be used to define the H\"older--Zygmund classes for all $\a>0$. By \rf{pbe}, for $\a>0$, the space
$\L_\a\big(\R^n\big)$ consists
of functions $f\in C\big(\R^n\big)$ such that
$$
|(\D_h^mf)(x)|\le\const|h|^\a,\quad x,~h\in\R^n,
$$
where $m$ is the smallest integer greater than $\a$.

We refer the reader to \cite{Pee} and \cite{T} for more detailed information on Besov spaces.

\medskip

{\bf 2.1. Operator ideals.}
In this section we give a brief introduction to quasinormed ideals of operators on Hilbert space. We refer the reader to \cite{AP3} for more detailed information.

Recall that a {\it quasinorm} $\|\cdot\|$ on a vector space
$X$ is an $\R_+$-valued function on $X$ such that $\|x\|>0$ unless $x=\0$;
 $\|\a x\|=|\a|\cdot\|x\|$,  $x\in X$, $\a\in\C$; and
 $\|x+y\|\le c\big(\|x\|+\|y\|)$, $x,\,y\in X$ for some $c>0$.

We say that a sequence $\{x_n\}$ of vectors of a {\it quasinormed space} $X$ converges to
$x\in X$ if $\lim\limits_{n\to\be}\|x_n-x\|=0$. A complete quasinormed space is called a {\it quasi-Banach space}.

If $X_1$ and $X_2$ are quasinormed spaces and
$T:X_1\to X_2$ is a continuous linear operator then the quasinorm of $T$ is, by definition, $\|T\|\df\sup\{\|Tx\|_{X_2}:~x\in X_1,~\|x\|_{X_1}\le1\}$.

For a bounded linear operator $T$ on Hilbert space, we consider its singular
values $s_j(T)$, $j\ge0$, 
$$
s_j(T)\df\inf\big\{\|T-R\|:~\rank R\le j\big\}.
$$
We also introduce the sequence $\{\s_n(T)\}_{n\ge0}$
defined by
\bay
\label{sin}
\s_n(T)\df\frac1{n+1}\sum_{j=0}^ns_j(T).
\ey
Note that unlike the singular values, the functional $\s_n(\cdot)$ are seminorms:
$$
\s_n(T+R)\le\s_n(T)+\s_n(R),\quad n\ge0,
$$
for arbitrary bounded linear operators $T$ and $R$. This follows easily from the equality 
$$
\sum_{j=0}^ns_j(T)=\sup\big|\trace(PUTQ)\big|,
$$
where the supremum is take over all rank $n$ projections $P$ and $Q$ and unitary operators $U$, see \cite{BS0}.

\medskip

{\bf Definition.} Let $\h$ be a Hilbert space and let ${\frak I}$ be a linear subset in the set $\B=\B(\h)$ of bounded linear operators on $\h$ that is equipped with a quasinorm $\|\cdot\|_{\frak I}$ that makes
$\fI$ a quasi-Banach space.
We say that ${\frak I}$ is a {\it quasinormed ideal} if for every $A$ and $B$ in $\B(\h)$ and for every
$T\in{\frak I}$,
\bay
\label{qni}
ATB\in{\frak I}\qm\|ATB\|_\fI\le\|A\|\cdot\|B\|\cdot\|T\|_\fI.
\ey
A quasinormed ideal $\fI$ is called a {\it normed ideal} if $\|\cdot\|_\fI$ is a norm.

Note that {\it we do not require that} $\fI\ne\B(\h)$.

\medskip

There exists an $[0,\be]$-valued function
$\Psi=\Psi_\fI$ on the set of nonincreasing sequences of nonnegative real numbers 
such that $T\in\fI$ if and only if $\Psi\big(s_0(T),s_1(T),s_2(T),\cdots~)<\be$ and
$$
\|T\|_\fI=\Psi\big(s_0(T),s_1(T),s_2(T),\cdots~),\quad T\in\fI,
$$
see \cite{GK}.

If $T$ is an operator from a Hilbert space $\h_1$ to a Hilbert space $\h_2$, we say that $T$ belongs to
$\fI$ if $\Psi\big(s_0(T),s_1(T),s_2(T),\cdots~)<\be$.

%For a quasinormed ideal $\fI$ and a positive number $p$, we define the quasinormed ideal $\fI^{\{p\}}$ by
%$$
%\fI^{\{p\}}=\left\{T:~\big(T^*T\big)^{p/2}\in\fI\right\},
%\quad\|T\|_{\fI^{\{p\}}}\df\left\|(T^*T\big)^{p/2}\right\|_\fI^{1/p}.
%$$

If $T$ is an operator on a Hilbert space $\h$ and $d$ is a positive integer, we denote by $[T]_d$
the operator $\bigoplus\limits_{j=1}^dT_j$ on the orthogonal sum $\bigoplus\limits_{j=1}^d\h$ of $d$ copies of $\h$, where
$T_j=T$, $1\le j\le d$. It is easy to see that
$$
s_n\big([T]_d\big)=s_{[n/d]}(T),\quad n\ge0,
$$
where $[x]$ denotes the largest integer that is less than or equal to $x$.

We denote by $\b_{\fI,d}$ the quasinorm of the transformer $T\mapsto[T]_d$
on $\fI$, i.e.,  
$$
\b_{\fI,d}=\sup\Big\{\big\|[T]_d\|_\fI:~\|T\|_\fI\le1\Big\},
$$
where the supremum is take over all operators $T$ in $\fI$ on Hilbert spaces.

Clearly, the sequence $\{\b_{\fI,d}\}_{d\ge1}$ is nondecreasing and {\it submultiplicative}, i.e.,
$$
\b_{\fI,d_1d_2}\le\b_{\fI,d_1}\b_{\fI,d_2}.
$$ 
It is well known (see e.g., \S\,3 of \cite{AP3}) that the last inequality implies that
$$
\lim_{d\to\be}\frac{\log\b_{\fI,d}}{\log d}=\inf_{d\ge2}\frac{\log\b_{\fI,d}}{\log d}.
$$

\medskip

{\bf Definition.} If $\fI$ is a quasinormed ideal, the number
$$
\b_\fI\df\lim_{d\to\be}\frac{\log\b_{\fI,d}}{\log d}=\inf_{d\ge2}\frac{\log\b_{\fI,d}}{\log d}
$$
is called the {\it upper Boyd index of} $\fI$.

\medskip

It is easy to see that $\b_\fI\le1$ for an arbitrary normed ideal $\fI$. It is also clear that $\b_\fI<1$
if and only if $\lim\limits_{d\to\be}d^{-1}\b_{\fI,d}=0$.

Note that the upper Boyd index does not change if we replace the initial quasinorm in the quasinormed ideal with an equivalent one that also satisfies \rf{qni}. %It is also easy to see that
%$$
%\b_{\fI^{\{p\}}}=p^{-1}\b_\fI.
%$$

The proof of the following fact can be found in \cite{AP3}, \S\,3.

\medskip

{\bf Theorem on ideals with upper Boyd index less than 1.}
{\em Let $\fI$ be a quasinormed ideal. The following are equivalent:

{\em(i)} $\b_\fI<1$;

{\em(ii)} for every nonincreasing sequence $\{s_n\}_{\ge0}$
of nonnegative numbers,
\bay
\label{xi}
\Psi_\fI\Big(\{\s_n\}_{n\ge0}\Big)\le\const\Psi_\fI\Big(\{s_n\}_{n\ge0}\Big),
\ey
where
$\s_n\df(1+n)^{-1}\sum\limits_{j=0}^ns_j$}.

\medskip

{\it We denote by}
$\bs{C}_\fI$ be the best possible constant in inequality \rf{xi}.
%Then (see \cite{AP3}, \S\,3)
%\bay
%\label{Ce}
%\bs{C}_\fI\le3\sum_{k=0}^\be2^{-k}\b_{\fI,2^k}.
%\ey
%
%\medskip

Let $\bS_p$, $0<p<\be$, be the Schatten--von Neumann class of operators $T$ on Hilbert space
such that
$$
\|T\|_{\bS_p}\df\left(\sum_{j\ge0}\big(s_j(T)\big)^p\right)^{1/p}<\be.
$$
This is a normed ideal for $p\ge1$. We denote by $\bS_{p,\be}$, $0<p<\be$, the ideal that consists of operators $T$ on Hilbert space such that
$$
\|T\|_{\bS_{p,\be}}\df\left(\sup_{j\ge0}(1+j)\big(s_j(T)\big)^p\right)^{1/p}.
$$
The quasinorm $\|\cdot\|_{p,\be}$ is not a norm, but it is equivalent to a norm if $p>1$.
It is easy to see that
$$
\b_{\bS_p}=\b_{\bS_{p,\be}}=\frac1p\,\,,\quad 0<p<\be.
$$
Thus $\bS_p$ and $\bS_{p,\be}$ with $p>1$ satisfy the hypotheses of the Theorem on ideals with upper Boyd index less than 1.

%It follows easily from \rf{Ce} that for $p>1$,
%$$
%\bs{C}_{\bS_p}\le3\big(1-2^{1/p-1}\big)^{-1}.
%$$

%Suppose now that $\fI$ is a quasinormed ideal of operators on Hilbert space.
%With a nonnegative integer $l$ we associate the ideal $^{(l)}\fI$
%that consists of all bounded linear operators on Hilbert space and
%is equipped with the norm
%$$
%\Psi_{^{(l)}\fI}(s_0,s_1,s_2,\cdots)=\Psi(s_0,s_1,\cdots,s_l,0,0,\cdots).
%$$
%It is easy to see that for every bounded operator $T$,
%\begin{align*}
%\|T\|_{^{(l)}\fI}&=\sup\big\{\|RT\|_\fI:~\|R\|\le1,~\rank R\le l+1\big\}
%\\[.2cm]
%&=\sup\big\{\|TR\|_\fI:~\|R\|\le1,~\rank R\le l+1\big\}.
%\end{align*}
%
%It is easy to verify (see \cite{AP3}, \S\,3) that if
% $\fI$ is a quasinormed ideal, then for all $l\ge0$,
%\bay
%\label{ClI}
%\bs{C}_{^{(l)}\fI}\le\bs{C}_\fI.
%\ey

%Note that if $\fI=\bS_p$, $p\ge1$, then $\bS_p^l\df{^{(l)}\bS_p}$
%is the normed ideal that consists of all bounded linear operators equipped with the norm
%$$
%\|T\|_{\bS_p^l}\df\left(\sum_{j=0}^l\big(s_j(T)\big)^p\right)^{1/p}.
%$$
%It is well known that $\|\cdot\|_{\bS_p^l}$ is a norm for $p\ge1$ (see \cite{BS0}).
%
%It is also well known (see \cite{AP3}, \S\,3) that
%\bay
%\label{rpq}
%\|T_1T_2\|_{\bS_r^l}\le\|T_1\|_{\bS_p^l}\|T_2\|_{\bS_q^l},
%\ey
%where $T_1$ and $T_2$ bounded operator on Hilbert space and $1/p+1/q=1/r$.

We say that a quasinormed ideal $\fI$ has {\it majorization property}  (respectively {\it weak majorization property}) if the conditions
$$
T_1\in\fI,\quad T_2\in\B,\quad
\mbox{and}\quad
\s_l(T_2)\le\s_l(T_1)\quad
\mbox{for all}\quad
l\ge0
$$
imply that
$$
T_2\in\fI\quad\mbox{and}\quad\|T_2\|_{\fI}\le\|T_1\|_{\fI}\quad(\text{respectively}\quad
\|T_2\|_{\fI}\le\const\|T_1\|_{\fI})
$$
(see \cite{GK}).
Note that if a quasinormed ideal $\fI$ has weak majorization property, then we can introduce on it the following new equivalent quasinorm:
$$
\|T\|_{\widetilde\fI}\df\sup\{\|R\|_{\fI}:~\s_l(R)\le\s_l(T)\,\,\,\text{for all}\,\,\,l\ge0\}
$$
such that
$(\fI,\|\cdot\|_{\widetilde\fI})$ has majorization property.

Every separable normed ideal and every normed ideal
that is dual to a separable normed ideal has majorization property, see \cite{GK}.
Clearly, $\bS_1\subset\fI$ for every quasinormed ideal $\fI$
with majorization property.
Note also that every quasinormed ideal $\fI$ with $\b_{\fI}<1$ has weak majorization property (see, for example, \S\,3 of \cite{AP3} and
\S\,3 of \cite{AP4}).

The following fact on interpolation properties of quasinormed ideals that have majorization property is well-known, see e.g., \cite{AP4}.

\medskip

{\bf Theorem on interpolation of quasinormed ideals.}
{\it Let $\fI$ be a quasinormed ideal with majorization property
and let $\frak A:\frak L\to\frak L$ be a linear transformer
on a linear subset $\frak L$ of $\B$ such that $\frak L\cap\bS_1$
is dense in $\bS_1$. Suppose that
$\|\frak A T\|\le\|T\|$ and $\|\frak A T\|_{\bS_1}\le\|T\|_{\bS_1}$
for all $T\in\frak L$. Then $\|\frak A T\|_{\fI}\le\|T\|_{\fI}$
for every  $T\in\frak L$}.

We refer the reader to \cite{GK} and \cite{BS0} for further information on singular values and normed ideals of operators on Hilbert space.

\medskip

{\bf 2.3. Double operator integrals.}
In this subsection we give a brief introduction to double  operator integrals. Double operator integrals appeared in the paper \cite{DK} by Daletskii and S.G. Krein. However, the beautiful theory of double operator integrals was developed later by Birman and Solomyak in \cite{BS1}, \cite{BS2}, and \cite{BS3}, see also their survey \cite{BS5}.

Let $(\X,E_1)$ and $(\Y,E_2)$ be spaces with spectral measures $E_1$ and $E_2$
on a Hilbert space $\h$. The idea of Birman and Solomyak is to define first
double operator integrals
\bay
\label{doi}
\int\limits_\X\int\limits_\Y\Phi(x,y)\,d E_1(x)T\,dE_2(y),
\ey
for bounded measurable functions $\Phi$ and operators $T$
of Hilbert Schmidt class $\bS_2$. Consider the spectral measure $\E$ whose values are orthogonal
projections on the Hilbert space $\bS_2$, which is defined by
$$
\E(\L\times\D)T=E_1(\L)TE_2(\D),\quad T\in\bS_2,
$$
$\L$ and $\D$ being measurable subsets of $\X$ and $\Y$. It was shown in \cite{BS} that $\E$ extends to a spectral measure on
$\X\times\Y$. If $\Phi$ is a bounded measurable function on $\X\times\Y$, we define the double operator integral \rf{doi} by
$$
\int\limits_\X\int\limits_\Y\Phi(x,y)\,d E_1(x)T\,dE_2(y)\df
\left(\,\,\int\limits_{\X\times\Y}\Phi\,d\E\right)T.
$$
Clearly,
$$
\left\|\int\limits_\X\int\limits_\Y\Phi(x,y)\,dE_1(x)T\,dE_2(y)\right\|_{\bS_2}
\le\|\Phi\|_{L^\be}\|T\|_{\bS_2}.
$$
If
$$
\int\limits_\X\int\limits_\Y\Phi(x,y)\,d E_1(x)T\,dE_2(y)\in\bS_1
$$
for every $T\in\bS_1$, we say that $\Phi$ is a {\it Schur multiplier of $\bS_1$ associated with
the spectral measures $E_1$ and $E_2$}.

In this case the transformer
\bay
\label{tra}
T\mapsto\int\limits_{\Y}\int\limits_{\X}\Phi(x,y)\,d E_2(y)\,T\,dE_1(x),\quad T\in \bS_2,
\ey
extends by duality to a bounded linear transformer on the space of bounded linear operators on $\h$
and we say that the function $\Psi$ on $\Y\times\X$ defined by
$$
\Psi(y,x)=\Phi(x,y)
$$
is {\it a Schur multiplier (with respect to $E_2$ and $E_1$) of the space of bounded linear operators}.
We denote the space of such Schur multipliers by $\fM(E_2,E_1)$.
The norm of $\Psi$ in $\fM(E_2,E_1)$ is, by definition, the norm of the
transformer \rf{tra} on the space of bounded linear operators.

In \cite{BS3} it was shown that if $A$ and $B$ are self-adjoint operators (not necessarily bounded) such that $A-B$ is bounded
 and if $f$ is a continuously differentiable
function on $\R$ such that the divided difference $\dg f$,
$$
\big(\dg f\big)(x,y)=\frac{f(x)-f(y)}{x-y},
$$
is a Schur multiplier
of $\bS_1$ with respect to the spectral measures of $A$ and $B$, then
$$
f(A)-f(B)=\iint\big(\dg f\big)(x,y)\,dE_{A}(x)(A-B)\,dE_B(y)
$$
and
$$
\|f(A)-f(B)\|\le\const\|f\|_{\fM(E_A,E_{B})}\|A-B\|,
$$
i.e., {\it $f$ is an operator Lipschitz function}.

There are different descriptions of the space $\fM(E_1,E_2)$ of Schur multipliers, see \cite{Pe2}. In particular, it follows from those descriptions that $\Phi\in\fM(E_1,E_2)$
if and only id $\Phi$ is a Schur multiplier of $\bS_1$ associated with $E_1$ and $E_2$.

In this paper we need the following easily verifiable sufficient condition:

\medskip

{\it If a function $\Phi$ on $\X\times\Y$ belongs to the {\it projective tensor
product}
$L^\be(E_1)\hat\otimes L^\be(E_2)$ of $L^\be(E_1)$ and $L^\be(E_2)$ (i.e., $\Phi$ admits a representation
$$
\Phi(x,y)=\sum_{j\ge0}\f_j(x)\psi_j(y),
$$
where $\f_j\in L^\be(E_1)$, $\psi_j\in L^\be(E_2)$, and
$$
\sum_{j\ge0}\|\f_j\|_{L^\be}\|\psi_j\|_{L^\be}<\be),
$$
then $\Phi\in\fM(E_1,E_2)$ and}
\bay
\label{dous}
\|\Phi\|_{\fM(E_1,E_2)}\le\sum_{j\ge0}\|\f_j\|_{L^\be}\|\psi_j\|_{L^\be}.
\ey

\medskip

For such functions $\Phi$ we have
$$
\int\limits_\X\int\limits_\Y\Phi(x,y)\,dE_1(x)T\,dE_2(y)=
\sum_{j\ge0}\left(\,\int\limits_\X\f_j\,dE_1\right)T\left(\,\int\limits_\Y\psi_j\,dE_2\right).
$$

It follows from
the Theorem on interpolation of quasinormed ideals (see Subsection 2.2) that if $\Phi\in\fM(E_1,E_2)$ and $\fI$ is a quasinormed ideal with majorization property, then
$$
T\in\fI\quad\Longrightarrow\quad
\int\limits_{\X}\int\limits_{\Y}\Phi(x,y)\,d E_1(x)\,T\,dE_2(y)\in\fI
$$
and
\bay
\label{intn}
\left\|\int\limits_{\X}\int\limits_{\Y}\Phi(x,y)\,d E_1(x)\,T\,dE_2(y)\right\|_\fI
\le\|\Phi\|_{\fM(E_1,E_2)}\|T\|_\fI.
\ey

Suppose now that $\D_1$ and $\D_2$ are Borel subsets of Euclidean spaces. We denote by $\fM_{\D_1,\D_2}$ the class of Borel functions $\Phi$ on $\D_1\times\D_2$ that belong to the space of Schur multipliers $\fM(E_1,E_2)$ for all Borel Spectral measures $E_1$ on $\D_1$ and $E_2$ on $\D_2$.  We put
$$
\|\Phi\|_{\fM_{\D_1,\D_2}}\df\sup\|\Phi\|_{\fM(E_1,E_2)},
$$
the supremum being taken over all such spectral measures $E_1$ and $E_2$. It follows from Theorem 2.2 of \cite{AP5} that the supremum is finite. 

The following result is well known and can be proved elementarily.

\begin{lem}
\label{bldia}
Let $\{\mQ_j\}$ be a family of disjoint Borel subsets of $\R^n$ and let $\{\mR_j\}$ be also a family of disjoint Borel subsets of $\R^m$. Suppose that $\Psi=\sum_j\Psi_j$, where $\Psi_j$ is a Borel function on $\R^n\times\R^m$ that is concentrated on 
$\mQ_j\times\mR_j$ (i.e., $\Psi_j\big|\R^n\times\R^m\setminus\mQ_j\times\mR_j=\0$) and
$\sup\|\Psi_j\|_{\fM_{\mQ_j,\mR_j}}<\be$.
Then $\Psi\in\fM_{\R^n,\R^m}$ and
$$
\|\Psi\|_{\fM_{\R^n\!,\R^m}}\le\sup_j\|\Psi_j\|_{\fM_{\mQ_j,\mR_j}}.
$$
\end{lem}

We need an elementary lemma that gives a sufficient condition for a function on the product of two cubes to be a Schur multiplier. To state the lemma, we introduce a piece of notation that will also be used in \S\,\ref{intf}. 

\medskip

{\bf Definition.} Let $Q$ be a cube in the Euclidean space $\R^d$ with center $c\in\R^d$ and let $K$ be a positive integer. We use the notation $K[Q]$ for the cube homothetic to $Q$ with respect to homothetic center $c$ and ratio $K$, i.e.,
\bay
\label{KQ}
K[Q]\df\big\{x\in\R^d:~c+K^{-1}(x-c)\in Q\big\}.
\ey

\begin{lem}
\label{osle}
Let $R=P\times Q\subset \R^{2n}=\R^n\times \R^n$ be a cube with sidelength $L$. 
Suppose that $\Psi$ is a $C^\infty$ function on $\frac32[R]$. For a multi-index $\a$ we put 
$$
C_\alpha=L^{|\alpha|}\max_{a\in\frac32[R]}\big|(D^\a\Psi)(a)\big|. 
$$
Then $\Psi\big|R\in\fM_{P,Q}$ and 
$$
\|\Psi\|_{\fM_{P,Q}}\le\const\max\Big\{\big|C_\a\big|:~|\a|\le2n+2\Big\}.
$$
\end{lem}

\Pf By applying a translation and a dilation, we may assume that $P=Q=\left[-\frac23\pi,\frac23\pi\right]$.

Let $\vk$ be a nonnegative $C^\be$ function on $\R^{2n}$ such that 
$$
\vk(a)=1\quad\mbox{for}\quad a\in\left[-\frac23\pi,\frac23\pi\right]^{2n}
\quad\mbox{and}\quad
\vk(a)=0\quad\mbox{for}\quad a\not\in\left[-\pi,\pi\right]^{2n}.
$$
Put $\Psi_\flat=\vk\Psi$. 

Consider the Fourier expansion of $\Psi_\flat$ on $[-\pi,\pi]^n\times[-\pi,\pi]^n$:
$$
\Psi_\flat(x,y)=
\sum_{r,s\in\Z^n}\widehat\Psi_\flat(r,s)
e^{{\rm i}\langle x,r\rangle}
e^{{\rm i}\langle y,s\rangle}.
$$
It follows from \rf{dous} that 
$$
\|\Psi\|_{\fM_{P,Q}}=\|\Psi_\flat\|_{\fM_{P,Q}}\le
\sum_{r,s\in\Z^n}\left|\widehat\Psi_\flat(r,s)\right|.
$$
The result is a consequence of the following elementary estimate of the moduli of Fourier coefficients:
$$
\left|\widehat\Psi_\flat(r,s)\right|
\le\const\max\Big\{\big|C_\a\big|:~|\a|\le2n+2\Big\}(|r|+|s|)^{-2n-2}.\quad\bl
$$

\

\section{\bf The techniques used in the case of normal operators\\ do not work when $\bs{n\ge3}$}
\setcounter{equation}{0}
\label{n=3}

\

In this section we are going review the approach of \cite{APPS2} for functions of normal operators and we will see that that approach does not generalize to the case of $n$-tuples of commuting self-adjoint operators with $n\ge3$. 

The proofs of the results of \cite{APPS2} for normal operators are based on the following formula: 
\begin{align}
\label{normal}
f(N_1)-f(N_2)=&\iint(\dg_xf)(z_1,z_2)\,dE_1(z_1)(A_1-A_2)\,dE_2(z_2)
\nonumber\\[.2cm]
&+
\iint(\dg_yf)(z_1,z_2)\,dE_1(z_1)(B_1-B_2)\,dE_2(z_2).
\end{align}
Here $N_1$ and $N_2$ are normal operators with bounded difference $N_1-N_2$,  
$A_j=\re N_j$, $B_j=\im N_j$, $x_j=\re z_j$, $y_j=\im z_j$, 
$E_j$ is the spectral measure of $N_j$,
$f$ is a bounded function on 
$\C=\R^2$ whose Fourier transform has compact support, 
%and the functions $\dg_xf$ and $\dg_yf$ on $\R^2\times\R^2$ are defined by
$$
(\dg_xf)(z_1,z_2)=\frac{f(x_1,y_1)-f(x_2,y_1)}{x_1-x_2},\quad z_1,~z_2\in\C,
$$
and
$$
(\dg_yf)(z_1,z_2)=\frac{f(x_2,y_1)-f(x_2,y_2)}{y_1-y_2},\quad z_1,~z_2\in\C.
$$
It was shown in \cite{APPS2} that under the above assumptions 
%under these assumptions the functions 
$\dg_xf$ and $\dg_yf$ belong to the space of Schur multipliers $\fM_{\R^2\!,\R^2}$ and formula \rf{normal} holds. This, in turn, was used in \cite{APPS2} to show that such functions $f$ are operator Lipschitz, i.e.,
$$
\|f(N_1)-f(N_2)\|\le\const\|N_1-N_2\|
$$
for arbitrary normal operators $N_1$ and $N_2$ such that $N_1-N_2$ is bounded.

However, 
%it turns out that 
in the case $n\ge3$ the situation is more complicated. Let $(A_1,A_2,A_3)$ and $(B_1,B_2,B_3)$ be triples of commuting self-adjoint operators. 
%We can assume for simplicity, that the operators are bounded. 
Suppose that $f$ is a bounded function on $\R^3$ whose Fourier transform has compact support. An analog of \rf{normal} for triples of commuting self-adjoint operators would be the formula
\begin{align}
\label{troika}
f(A_1,A_2,A_3)-f(B_1,B_2,B_3)=&\iint(\dg_1f)(x,y)\,dE_A(x)(A_1-B_1)\,dE_B(y)\nonumber\\[.2cm]
&+\iint (\dg_2f)(x,y)\,dE_A(x)(A_2-B_2)\,dE_B(y)\nonumber\\[.2cm]
&+\iint (\dg_3f)(x,y)\,dE_A(x)(A_3-B_3)\,dE_B(y),
\end{align}
where
$$
(\dg_1f)(x,y)=\frac{f(x_1,x_2,x_3)-f(y_1,x_2,x_3)}{x_1-y_1},\quad
(\dg_2f)(x,y)=\frac{f(y_1,x_2,x_3)-f(y_1,y_2,x_3)}{x_2-y_2},
$$
$$
(\dg_3f)(x,y)=\frac{f(y_1,y_2,x_3)-f(y_1,y_2,y_3)}{x_3-y_3},
\quad x=(x_1,x_2,x_3),\quad y=(y_1,y_2,y_3),
$$
and $E_A$ and $E_B$ are the joint spectral measures of the triples $(A_1,A_2,A_3)$ and 
$(B_1,B_2,B_3)$ on the Euclidean space $\R^3$.

It can easily be shown that \rf{troika} holds if
the functions $\dg_1f$, $\dg_2f$, and $\dg_3f$ belong to the space of Schur multipliers $\fM_{\R^3\!,\R^3}$ which would imply that $f$ is an operator Lipschitz function.

The methods of \cite{APPS2} allow us to prove that $\dg_1f$ and $\dg_3f$ do 
belong to the space of Schur multipliers $\fM_{\R^3\!,\R^3}$. However, as the next result shows, the function $\dg_2f$ does not have to be in $\fM_{\R^3\!,\R^3}$, and so formula \rf{troika} cannot be used to prove that bounded functions on $\R^3$ with compactly supported Fourier transform must be operator Lipschitz.

\begin{thm}
\label{kontr}
Suppose that $g$ is a bounded continuous function on $\R$ such that the Fourier transform of $g$ has compact support and is not a measure. Let $f$ be the function on $\R^3$ defined by
\bay
\label{ef}
f(x_1,x_2,x_3)=g(x_1-x_3)\sin x_2,\quad x_1,\,x_2,\,x_3\in\R.
\ey 
Then $f$ is a bounded function on $\R^3$ whose Fourier transform has compact support, but $\dg_2f\not\in\fM_{\R^3\!,\R^3}$.
\end{thm}

To construct  a function $g$ satisfying the hypothesis of Theorem \ref{kontr}, one can take, for example, the function $g$ defined by
$$
g(x)=\int_0^xt^{-1}\sin t\,dt,\quad x\in\R.
$$
Obviously, $g$ is bounded 
and
%$$
%\lim_{t\to-\be}g(t)\ne\lim_{t\to\be}g(t).
%$$
and its Fourier transform $\F g$ satisfies the equality:
$$
(\F g)(t)=\frac ct
$$
for a nonzero constant $c$and  sufficiently small positive $t$. Clearly, this implies that $\F g$ is not a measure.

\medskip

\Pf Let $f$ be the function defined by \rf{ef}. Then $f$ is a bounded function whose Fourier transform is compactly supported. 

We have
$$
(\dg_2f)(x,y)=g(y_1-x_3)\frac{\sin x_2-\sin y_2}{x_2-y_2}, 
\quad x=(x_1,x_2,x_3),\quad y=(y_1,y_2,y_3).
$$
Since $g$ is not the Fourier transform of a measure, it is well known that the function
$$
(x,y)\mapsto g(x-y),\quad x,\,y\in\R,
$$
does not belong to the space $\fM_{\R,\R}$. Moreover, there exists a continuous function
$\vk$ on $\R$ such that $\vk$ is the Fourier transform of an $L^\be$ function (and so the operator of convolution with $\vk$ is bounded on $L^2(\R)$) such that  convolution with $g\vk$ is not a bounded linear operator on $L^2(\R)$. Let us fix such a function $\vk$.

To prove that $\dg_2f\not\in\fM_{\R^3\!,\R^3}$, consider kernels of the form
$$
k_n(x_1,x_2,x_3,y_1,y_2,y_3)=\vk(y_1-x_3)\,\chi_n(x_3,y_1)\,\xi(x_1,y_3)\,\eta(x_2,y_2),
$$
where $\chi_n$ is the characteristic function of $[-n,n]^2$, and $\xi$ and $\eta$ are  nonzero functions in $L^2(\R^2)$.

Clearly, the integral operators on $L^2(\R^3)$ with kernel functions $k_n$ are Hilbert Schmidt operators with uniformly bounded operator norms. On the other hand, it is also easy to verify that the operator norms of the integral operators on $L^2(\R^3)$ with kernel functions $(\dg_2f)k_n$ tend to infinity. This implies the desired conclusion. $\bl$

\

\section{\bf An integral formula}
\setcounter{equation}{0}
\label{intf}

\

We have seen in the previous section that for $n>2$, the representation
$$
f(x_1,\cdots,x_n)-f(y_1,\cdots,y_n)=\sum_{j=1}^n(x_j-y_j)(\dg_jf)(x_1,\cdots,x_n,y_1,\cdots,y_n)
$$
(the divided differences $\dg_jf$ can be defined by analogy with the definition given in \S\,\ref{n=3} for $n=3$)
does not help to establish that an $L^\be$ function $f$ on $\R^n$ whose Fourier transform has compact support is operator Lipschitz. The reason is that if $2\le j \le n-1$, the function 
$\dg_jf$ does not have to be a Schur multiplier. 

In this section we show that we can successfully replace the divided differences $\dg_jf$
with other functions that are Schur multipliers. We find, for an
 $L^\be$ function $f$ on $\R^n$ whose Fourier transform has compact support, a representation of the form
\bay
\label{pred} 
f(x_1,\cdots,x_n)-f(y_1,\cdots,y_n)=\sum_{j=1}^n(x_j-y_j)\Psi_j(x_1,\cdots,x_n,y_1,\cdots,y_n),
\ey
where the $\Psi_j$ are Schur multipliers.

This allows us to represent the difference $f(A_1,\cdots,A_n)-f(B_1,\cdots,B_n)$ for commuting $n$-tuples $(A_1,\cdots,A_n)$ and $(B_1,\cdots,B_n)$ in terms of the sum of certain double operator integrals, which will be used in the next section to obtain operator Lipschitz estimates.

\begin{thm}
\label{fsi}
Let $\s>0$ and let $f$ be a function in $L^\be(\R^n)$ whose Fourier transform is supported in the ball
$\{\xi\in\R^n:~\|\xi\|\le\s\}$. Then there exist functions $\Psi_j$, $1\le j\le n$, on $\R^n\times\R^n$
such that $\Psi_j\in\fM_{\R^n\!,\R^n}$, representation {\em\rf{pred}} holds,
and
$$
\|\Psi_j\|_{\fM_{\R^n\!,\R^n}}\le C_n\s\|f\|_{L^\be(\R^n)}
$$
for some positive number $C_n$ that does not depend on $\s$ and $f$.
\end{thm}

Note that by {\it Bernstein's inequality}, for all multi-indices $\a$,
\bay
\label{Ber}
\|D^\a f\|_{L^\be(\R^n)}\le\s^{|\a|}\|f\|_{L^\be(\R^n)}.
\ey

Recall that a {\it dyadic} interval is an interval of the form $\big[j2^k,(j+1)2^k\big)$, $j,\,k\in\Z$.
By a {\it dyadic cube} in $\R^d$ we mean a cube $I_1\times\cdots\times I_d$ such that each $I_j$ is a dyadic interval.

If $\mC$ is a dyadic cube, there is a unique dyadic cube whose sidelength is twice as large as the sidelength of $\mC$. We call it the {\it parent cube of} $\mC$.

We say that a dyadic cube $\mC$ in $\R^{2n}=\R^n\times\R^n$ is {\it admissible} if either its sidelength $l(\mC)$ is equal to $1$ or $l(\mC)>1$ and the interior of the cube $2[\mC]$ does not intersect the diagonal $\{(x,x):~x\in\R^n\}$ (see
\rf{KQ} for the definition of $2[\mC]$). An admissible cube is called {\it maximal} if it is not a proper subset of another admissible cube.

Consider all maximal admissible cubes. It is easy to see that they are disjoint and cover $\R^{2n}$. 

We need the following fact:

\begin{lem}
\label{kuby}
If $\mQ$ is a dyadic cube in $\R^n$, then there can be at most
$6^n$ dyadic cubes $\mR$ in $\R^n$ such that $\mQ\times\mR$ is a maximal admissible cube.
Similarly, for each dyadic cube $\mR$ in $\R^n$ there can be at most $6^n$ dyadic cube $\mQ$ in $\R^n$ such that $\mQ\times\mR$ is a maximal admissible cube. 
\end{lem}

\newcommand{\wh}{\widehat}

\Pf Let $\mC=\mQ\times \mR$ be a maximal admissible cube of sidelength $l(\mC)=l(\mQ)=l(\mR)=k$. Denote by $\wh\mQ$ and $\wh\mR$ the parent cubes of 
$\mQ$ and $\mR$. Then
the cube $\wh\mC=\wh\mQ\times \wh\mR$ is not admissible. This means that $(x,x)\in2[\wh\mC]$ for some $x\in\R^n$, i.e., that
$2[\wh\mQ]\cap 2[\wh\mR]\ne\varnothing$. Since  $l(\wh\mQ)=l(\wh\mR)=2k$, the last condition is 
equivalent to $3[\wh\mQ]\cap [\wh\mR]\ne\varnothing$. Note now that $3[\wh\mQ]$ consists of $3^n$
dyadic cubes of sidelength $2k$ and that two different dyadic cubes of the same size are disjoint.
Thus, for every given dyadic cube $\mQ$, there can be at most $3^n$ options for $\wh\mR$ and, thereby,
at most $6^n$ options for $\mR$. $\bl$

\medskip

For $k=1,2,4,8,\cdots$, we denote by $\cd_k$ the set of maximal dyadic cubes of sidelength $k$.

\medskip

{\bf Remark.} It follows now from Lemmata  \ref{kuby} and \ref{bldia} that if $k=2^m$, $m\in\Z_+$, and for $\mC\in\cd_k$, $\Psi_\mC$ is a Borel function vanishing outside $\mC$, then  
$$
\left\|\sum_{\mC\in\cd_k}\Psi_\mC\right\|_{\fM_{\R^n\!,\R^n}}
\le6^{2n}\sup\left\{\big\|\Psi_\mC\big\|_{\fM_{\mQ,\mR}}:~\mC\in\cd_k\right\}.
$$

\medskip 

{\bf Proof of Theorem \ref{fsi}.} It is easy to see that by rescaling, we may assume that $\s=1$.

We construct functions $\Psi_j^{[k]}$, $1\le j\le n$, $k=2^m$, $m\in\Z_+$, such that $\Psi_j^{[k]}$ is concentrated on $\bigcup\limits_{\mC\in\cd_k}\mC$, 
$$
\sum_{m\ge0}\Big\|\Psi^{[2^m]}_j\Big\|_{\fM_{\R^n\!,\R^n}}<\be,\quad1\le j\le n,
$$
and
$$ 
f(x)-f(y)=\sum_{j=1}^n(x_j-y_j)\Psi^{[k]}_j(x,y),\quad\mbox{for}\quad(x,y)\in\bigcup_{\mC\in\cd_k}\mC,
$$
where $x=(x_1,\cdots,x_n)$ and $y=(y_1,\cdots,y_n)$. Then it remains to put
\bay
\label{psij}
\Psi_j=\sum_{m\ge0}\Psi^{[2^m]}_j.
\ey

Next, by the Remark following the proof of Lemma \ref{kuby}, to estimate $\Big\|\Psi^{[k]}_j\Big\|_{\fM_{\R^n\!,\R^n}}$, it suffices to estimate the Schur multiplier norm of $\Psi^{[k]}_j$
on each cube $\mC$ in $\cd_k$.

Consider first the case $k=1$. Let $\mC=\mQ\times\mR\in\cd_1$. We have
$$
f(x)-f(y)=\sum_{j=1}^n(x_j-y_j)\int_0^1 (D_jf)((1-t)x+ty)\,dt,
\quad x\in\mQ,~y\in\mR,
$$
where $D_jf$ is the $j$th partial derivative of $f$.
Put
$$
\Phi_j(x,y)=\int_0^1(D_jf)((1-t)x+ty)\,dt,\quad(x,y)\in\frac32[\mC].
$$

It follows from \rf{Ber} and from Lemma \ref{osle} that
$$
\|\Phi_j\|_{\fM_{\mQ,\mR}}\le\const\|f\|_{L^\be(\R^n)}.
$$
We can put now
$$
\Psi^{[1]}_j=\sum_{\mC\in\cd_1}\chi_{_\mC}\Phi_j.
$$
By Remark 2,
$$
\Big\|\Psi^{[1]}_j\Big\|_{\fM_{\R^n\!,\R^n}}\le\const\|f\|_{L^\be(\R^n)}.
$$

We proceed now to the case $k>1$. Suppose that $\mC=\mQ\times\mR\in\cd_k$. 
Let $\o$ be a $C^\be$ nonnegative even function on $\R$ such that 
$$
\o(t)=0,\quad t\in\left[-\frac12,\frac12\right],\qquad\mbox{and}\qquad\o(t)=1,\quad t\not\in[-1,1].
$$
For $x=(x_1,\cdots,x_n)$ and $y=(y_1,\cdots,y_n)$ in $\R^n$, we put
$$
\Phi_j(x,y)=\o\left(\frac{x_j-y_j}{k}\right)\quad
\mbox{and}\quad
\Phi=\sum_{j=1}^n\Phi_j,
$$
and define the functions $\Xi_j$, $1\le j\le n$, by
$$
\Xi_j(x,y)=\left\{\begin{array}{ll}
\frac{1}{x_j-y_j}\cdot\frac{\Phi_j(x,y)}{\Phi(x,y)},&x_j\ne y_j,\\[.2cm]
0,&x_j=y_j.\end{array}\right.
$$
It is easy to see that
$$
\sum_{j=1}^n(x_j-y_j)\Xi_j(x,y)=1,\quad (x,y)\in\frac32[\mC],
$$
and so
$$
f(x)-f(y)=\sum_{j=1}^n(x_j-y_j)\big(f(x)-f(y)\big)\Xi_j(x,y),\quad (x,y)\in\frac32[\mC].
$$
It is easy to verify that the function $k\Xi_j$ satisfies the hypotheses of Lemma \ref{osle}. Thus $\|\Xi_j\|_{\fM_{\mQ,\mR}}\le\const k^{-1}$. We define the function $\Psi_j^{[k]}$ on $\mC$ by
$$
\Psi_j^{[k]}(x,y)=\big(f(x)-f(y)\big)\Xi_j(x,y),\quad(x,y)\in\mC.
$$
Since $\|f\|_{L^\be(\R^n)}\le1$, it follows that
$$
\left\|\Psi_j^{[k]}\right\|_{\fM_{\mQ,\mR}}\le\const k^{-1},
$$
and by Remarks 1 and 2,
$$
\left\|\Psi_j^{[k]}\right\|_{\fM_{\R^n,\R^n}}\le\const\sup_{\mC=\mQ\times\mR\in\cd_k}
\left\|\Psi_j^{[k]}\right\|_{\fM_{\mQ,\mR}}\le\const k^{-1}.
$$
It remains to observe that the function $\Psi_j$ defined by \rf{psij} satisfies
$$
\|\Psi\|_{\fM_{\R^n,\R^n}}\le\const\sum_{m\ge0}\left\|\Psi_j^{[2^m]}\right\|_{\fM_{\R^n,\R^n}}\le\const\sum_{m\ge0}2^{-m}<\be.\quad\bl
$$

Theorem \ref{fsi} easily implies that $L^\be$ functions whose Fourier transform has compact support are operator Lipschitz. More precisely, the following result holds:

\begin{thm}
\label{OLcs}
Let $f$ satisfy the hypotheses of Theorem {\em\ref{fsi}}. Suppose that $A_1,\cdots,A_n$ and $B_1,\cdots,B_n$ are $n$-tuples of commuting self-adjoint operators on Hilbert space 
such that $A_j-B_j$ is bounded for $1\le j\le n$, and let $E_A$ and $E_B$ be their joint spectral measures on $\R^n$. Then
\begin{align}
\label{nintf}
f(A_1,\cdots,A_n)-f(B_1,\cdots,B_n)=\sum_{j=1}^n\iint\Psi_j(x,y)\,
dE_A(x)(A_j-B_j)\,dE_B(y),
\end{align}
where the functions $\Psi_j$ satisfy the conclusion of Theorem {\em\ref{fsi}}.
\end{thm}

\Pf Let us first establish formula \rf{nintf} for bounded operators $A_1,\cdots,A_n$ and $B_1,\cdots,B_n$. We have
$$
\iint\Psi_j(x,y)\,dE_A(x)(A_j-B_j)\,dE_B(y)=\iint\Psi_j(x,y)(x_j-y_j)\,dE_A(x)\,dE_B(y).
$$
Thus, using \rf{pred}, we obtain
\begin{align*}
\sum_{j=1}^n\iint\Psi_j(x,y)\,dE_A(x)(A_j-B_j)\,dE_B(y)&=
\sum_{j=1}^n\iint\Psi_j(x,y)(x_j-y_j)\,dE_A(x)\,dE_B(y)\\[.2cm]
&=\iint\big(f(x)-f(y)\big)\,dE_A(x)\,dE_B(y)\\[.2cm]
&=f(A_1,\cdots,A_n)-f(B_1,\cdots,B_n).
\end{align*}

In the general case when the operators are not necessarily bounded we argue as in the proof of Theorem 5.2 of \cite{APPS2}: for a positive integer $k$, we define the orthogonal projections
$$
P_k=E_A\left(\left\{x\in\R^n:~\sum_{j=1}^n|x_j|\le k\right\}\right)\quad\mbox{and}\quad
Q_k=E_B\left(\left\{x\in\R^n:~\sum_{j=1}^n|x_j|\le k\right\}\right).
$$
Applying equality\rf{nintf} for the $n$-tuples of bounded commuting self-adjoint operators
$(A_{1,k},\cdots,A_{n,k})$ and $(B_{1,k},\cdots,B_{n,k})$,
$$
A_{j,k}\df P_kA_j\quad\mbox{and}\quad B_{j,k}\df Q_kB_j,\quad 1\le j\le n,
$$
we obtain
\begin{align*}
P_k\big(f(A_1,\cdots,A_n)&-f(B_1,\cdots,B_n)\big)Q_k\\[.2cm]
&=P_k\big(P_kf(A_1,\cdots,A_n)-f(B_1,\cdots,B_n)Q_k\big)Q_k\\[.2cm]
&=\sum_{j=1}^nP_k\left(\iint\Psi_j(x,y)\,dE_A(x)(A_{j,k}-B_{j,k})\,dE_B(y)
\right)Q_k.
\end{align*}
To obtain formula \rf{nintf}, we take the limit as $k\to\be$ in the strong operator topology, see details in the proof of Theorem 5.2 of \cite{APPS2}.
$\bl$

\begin{thm}
\label{Besif}
Let $f\in B_{\be,1}^1(\R^n)$ and let $A_1,\cdots,A_n$ and $B_1,\cdots,B_n$ be $n$-tuples of commuting self-adjoint operators
such that the operators $A_j-B_j$, $1\le j\le n$, are bounded.
Then formula {\em\rf{nintf}} holds.
\end{thm}

\Pf Put $f_l=f*W_l$, $l\in\Z$. Clearly, $f_l$ satisfies the hypotheses of Theorem \ref{fsi}
with $\s=2^{l+1}$.
To prove formula \rf{nintf}, it suffices to apply formula \rf{nintf} to each function $f_l$ (this can be done in view of Theorem \ref{OLcs}) and take the sum over $l\in\Z$. $\bl$

\

\section{\bf Operator Lipschitzness}
\setcounter{equation}{0}
\label{Olss}

\

In this section we show that functions in the Besov class $B_{\be,1}^1(\R^n)$ are operator Lipschitz. Moreover, 
for functions in $B_{\be,1}^1(\R^n)$ we also obtain Lipschitz estimates in operator ideal norms under a very mild assumption on the ideal.

As in the case of functions on $\R$ and $\R^2$, a function $f$ on $\R^n$ is called {\it operator Lipschitz} if 
\bay
\label{oLn}
\big\|f(A_1,\cdots,A_n)-f(B_1,\cdots,B_n)\big\|
\le\const\max_{1\le j\le n}\|A_j-B_j\|,
\ey
whenever $(A_1,\cdots,A_n)$ and $(B_1,\cdots,B_n)$ are $n$-tuples of commuting self-adjoint operators such that the operators $A_j-B_j$ are bounded, 
$1\le j\le n$. The operator Lipschitz norm $\|f\|_{\rm OL}$ of $f$ is, by definition, the minimal possible constant in inequality \rf{oLn}.

\begin{thm}
\label{sigOL}
Let $f$ be an $L^\be$ function on $\R^n$ whose Fourier transform  
is supported in the ball $\{\xi\in\R^n:~\|\xi\|\le\s\}$.
Then $f$ is operator Lipschitz and
\bay
\label{OLsi}
\|f\|_{\rm OL}\le c_n\s\|f\|_{L^\be}
\ey
for some positive number $c_n$.
\end{thm}

\Pf Let $A_1,\cdots,A_n$ and $B_1,\cdots,B_n$ be $n$-tuples of commuting self-adjoint operators on Hilbert space such that $A_j-B_j$ is bounded for $1\le j\le n$.
It follows from formula \rf{nintf} and from Theorem \ref{fsi} that
\begin{align*}
\big\|f(A_1,\cdots,A_n)-f(B_1,\cdots,B_n)\big\|&\le\sum_{j=1}^n\|\Psi_j\|_{\fM_{\R^n,\R^n}}\|A_j-B_j\|\\[.2cm]
&\le\const\s\|f\|_{L^\be}\max_{1\le j\le n}\|A_j-B_j\|.
\end{align*}
This implies the result. $\bl$ 

\begin{thm}
\label{BesOL}
Let $f\in B_{\be,1}^1(\R^n)$. Then $f$ is operator Lipschitz and
$$
\|f\|_{\rm OL}\le c_n\|f\|_{B_{\be,1}^1}
$$
for some positive number $c_n$.
\end{thm}

\Pf Indeed, let $f_l=f*W_l$, $l\in\Z$. Suppose that $A_1,\cdots,A_n$ and $B_1,\cdots,B_n$ are $n$-tuples of commuting self-adjoint operators on Hilbert space such that $A_j-B_j$ is bounded for $1\le j\le n$. We have
\begin{align}
\label{opn}
\big\|f(A_1,\cdots,A_n)-f(B_1,\cdots,B_n)\big\|
&\le\sum_{l\in\Z}\big\|f_l(A_1,\cdots,A_n)-f_l(B_1,\cdots,B_n)\big\|\nonumber\\[.2cm]
&\le\const\sum_{l\in\Z}2^l\|f_l\|_{L^\be}\max_{1\le j\le n}\|A_j-B_j\|\nonumber\\[.2cm]
&\le\const\|f\|_{B_{\be,1}^1}\max_{1\le j\le n}\|A_j-B_j\|
\end{align}
by \rf{OLsi}. $\bl$

We proceed now to Lipschitz type estimates in operator ideal (quasi)norms.

\begin{thm}
\label{yadn}
Let $\fI$ be a quasinormed ideal of operators on Hilbert space that has majorization property.
Suppose that $A_1,\cdots,A_n$ and $B_1,\cdots,B_n$ are $n$-tuples of commuting self-adjoint operators on Hilbert space such that $A_j-B_j\in\fI$ for $1\le j\le n$.

If $f$ is an $L^\be$ function on $\R^n$ whose Fourier transform is supported in the ball
\lb$\{\xi\in\R^n:~\|\xi\|\le\s\}$, then $f(A_1,\cdots,A_n)-f(B_1,\cdots,B_n)\in\fI$ and
\bay
\label{siId}
\|f(A_1,\cdots,A_n)-f(B_1,\cdots,B_n)\|_\fI\le c_n\s\|f\|_{L^\be}
\max_{1\le j\le n}\|A_j-B_j\|_\fI
\ey
for some positive number $c_n$.

If $f$ is a function in the Besov class $B_{\be1}^1\big(\R^n\big)$,
then $f(A_1,\cdots,A_n)-f(B_1,\cdots,B_n)\in\fI$ and 
$$
\|f(A_1,\cdots,A_n)-f(B_1,\cdots,B_n)\|_\fI\le c_n\,\|f\|_{B_{\be1}^1}
\max_{1\le j\le n}\|A_j-B_j\|_\fI.
$$
for some positive number $c_n$.
\end{thm}

\Pf The result follows immediately from Theorems \ref{OLcs} and \ref{Besif},
and from \rf{intn} $\bl$

\medskip

In particular, we can apply Theorem \ref{yadn} for trace class perturbations.

\begin{cor}
\label{BesTr}
Let $f\in B_{\be,1}^1(\R^n)$ and suppose that $A_1,\cdots,A_n$ and $B_1,\cdots,B_n$ are $n$-tuples of commuting self-adjoint operators on Hilbert space such that $A_j-B_j\in\bS_1$  for $1\le j\le n$. Then $f(A_1,\cdots,A_n)-f(B_1,\cdots,B_n)\in\bS_1$ and
$$
\big\|f(A_1,\cdots,A_n)-f(B_1,\cdots,B_n)\big\|_{\bS_1}
\le\const\|f\|_{B_{\be,1}^1}\max_{1\le j\le n}\|A_j-B_j\|_{\bS_1}.
$$ 
\end{cor}

\Pf The proof is the same as the proof of \rf{opn}. $\bl$

\medskip

Note that Theorem \ref{yadn} also implies Lipschitz estimates in the Schatten--von Neumann norm $\bS_p$ with $p\ge1$. However, for $p\in(1,\be)$, a much stronger result was obtained in \cite{KPSS}. Namely, it was shown in \cite{KPSS} that if $1<p<\be$ and 
$A_1,\cdots,A_n$ and $B_1,\cdots,B_n$ are $n$-tuples of commuting self-adjoint operators on Hilbert space such that $A_j-B_j\in\bS_p$, then for every Lipschitz function $f$ on $\R^n$, $$
f(A_1,\cdots,A_n)-f(B_1,\cdots,B_n)\in\bS_p
$$
and
$$
\big\|f(A_1,\cdots,A_n)-f(B_1,\cdots,B_n)\big\|_{\bS_p}
\le\const\|f\|_{\rm Lip}\max_{1\le j\le n}\|A_j-B_j\|_{\bS_p}.
$$

Note also that in \cite{NP0} it is proved that for a Lipschitz function $f$ on $\R$ and for self-adjoint operators $A$ and $B$ with $\rank(A-B)<\be$, the operator $f(A)-f(B)$ belongs to the quasi-normed ideal $\bS_{1,\be}$, i.e.,
$$
s_m\big(f(A)-f(B)\big)\le\const(1+m)^{-1},\quad m\ge0.
$$
It is still unknown whether the same conclusion holds for Lipschitz functions $f$ under the assumption that $A-B\in\bS_1$.

\

\section{\bf H\"older classes and arbitrary moduli of continuity}
\setcounter{equation}{0}
\label{HcMc}

\

The purpose of this section is to obtain sharp estimates for the operator norms of $f(A_1,\cdots,A_n)-f(B_1,\cdots,B_n)$
for $n$-tuples $(A_1,\cdots,A_n)$ and $(B_1,\cdots,B_n)$ of commuting self-adjoint operators in the case when $f$ belongs to the 
H\"older class $\L_\a(\R^n)$ or, more generally, $f$ belongs to the class $\L_\o(\R^n)$, where $\o$ is an arbitrary modulus of continuity. 

We establish analogs of the results of \cite{AP2} for perturbations of functions of self-adjoint operators (this corresponds to the case $n=1$). Recall that similar results for perturbations of functions of normal operators were obtained in \cite{APPS2} (this corresponds to the case $n=2$). We generalize in this section the results of \cite{AP2} and
\cite{APPS2} to the case of arbitrary $n$.

The crucial step to obtain the above mentioned results of \cite{AP2} and
\cite{APPS2} was the fact that if $f$ is an $L^\be$ function on $\R$ (or on $\R^2$) whose Fourier transform has compact support, then $f$ is operator Lipschitz. In \S\,\ref{intf} of this paper we have shown that the same is true for $L^\be$ functions on $\R^n$, $n\ge1$, whose Fourier transform has compact support. Using this result, we can obtain the results of this section in the same way as in \cite{AP2} (or in \cite{APPS2}). 
Thus we only state the results and refer the reader to \cite{AP2} (or to \cite{APPS2}) to see how the results can be proved.

\begin{thm}
\label{0a1}
Let $n$ be a positive integer. There exists a positive number $c_n$ such that for every $\a\in(0,1)$ and for every $f\in\L_\a(\R^n)$,
\bay
\label{a^-1}
\|f(A_1,\cdots,A_n)-f(B_1,\cdots,B_n)\|\le c_n(1-\a)^{-1}\max_{1\le j\le n}\|A_j-B_j\|^\a,
\ey
whenever $(A_1,\cdots,A_n)$ and $(B_1,\cdots,B_n)$ are $n$-tuples of commuting self-adjoint operators.
\end{thm}

{\bf Remark.} Let $\a\in(0,1)$. Put
$$
{\frak h}_\a\df\sup\left\{\frac{\|f(A)-f(B)\|}{\|A-B\|^\a}
\right\},
$$
where the supremum is taken over all functions $f\in\L_\a(\R)$ with $\|f\|_{\L_\a}\le1$
and all bounded self-adjoint operators $A$ and $B$ such that $A\ne B$.
It was proved in Theorem 7.1 of \cite{AP6} that 
$$
{\frak h}_\a\ge c(1-\a)^{-1/2},\quad\a\in(0,1),
$$
for some positive number $c$.
It follows that in inequality \rf{a^-1} one cannot replace $(1-\a)^{-1}$ with anything better than $(1-\a)^{-1/2}$.

\medskip

We proceed now to the case of an arbitrary modulus of continuity $\o$. The function $\o_*$ on $(0,\be)$ is defined by
\bay
\label{o_*}
\o_*(\d)\df\d\int_\d^\be\frac{\o(t)}{t^2}\,dt,\quad\d>0.
\ey

\begin{thm}
\label{amC}
Let $n$ be a positive integer. There exists a positive number $c_n$ such that for every modulus of continuity $\o$ and for every $f\in\L_\o(\R^n)$,
$$
\|f(A_1,\cdots,A_n)-f(B_1,\cdots,B_n)\|\le c_n\|f\|_{\L_\o}\max_{1\le j\le n}\o_*\left(\|A_j-B_j\|\right),
$$
whenever $(A_1,\cdots,A_n)$ and $(B_1,\cdots,B_n)$ are $n$-tuples of commuting self-adjoint operators.
\end{thm}

In the case $\o(\t)=t$ (i.e., $\L_\o(\R^n)$ is the space ${\rm Lip}(\R^n)$ of Lipschitz functions), the function $\o_*$ identically takes value $\be$ on $(0,\be$), and so Theorem \ref{amC} does not give any estimate for Lipschitz functions. However, in the case when the joint spectra of $n$-tuples $(A_1,\cdots,A_n)$ and $(B_1,\cdots,B_n)$ are contained in a fixed bounded subset $K$ of $\R^n$ we can replace the initial modulus of continuity $\o(t)=t$ with the modulus of continuity defined by
$$
\o(\d)=\left\{\begin{array}{ll}\d,&\d\le d,\\[.2cm]d,&\d>d,
\end{array}\right.
$$
and apply Theorem \ref{amC}. This leads to the following result.

\begin{thm}
\label{Lippert}
Let $n$ be a positive integer. There exists a positive number $c_n$ such that  for every $f\in\Li(\R^n)$,
$$
\|f(A_1,\cdots,A_n)-f(B_1,\cdots,B_n)\|\le 
c_n\|f\|_{\Li}\max_{1\le j\le n}\|A_j-B_j\|\!
\left(\!1+\log\frac{d}{\!\max\limits_{1\le j\le n}\!\|A_j-B_j\|}\!\right),
$$
whenever $(A_1,\cdots,A_n)$ and $(B_1,\cdots,B_n)$ are $n$-tuples of commuting self-adjoint operators whose joint spectra are contained in a bounded set $K$ and $d$ is the diameter of $K$.
\end{thm}

The proof is similar to the proof of Theorem 8.4 of \cite{APPS2}, where the case of functions of normal operators is considered.

\

\section{\bf Estimates in ideal norms}
\setcounter{equation}{0}
\label{Idealc}

\

In this section we estimate ideal norms of perturbations of functions of $n$-tuples of commuting self-adjoint operators. Such results were obtained earlier in \cite{AP3} in the case of functions of self-adjoint operators and in \cite{APPS2} in the case of functions of normal operators. The results of this section can be deduced from inequality \rf{siId} in the same way as the results of \cite{AP3} for functions of self-adjoint operators and the results of \cite{APPS2} for functions of normal operators.

\begin{thm}
\label{Spnki}
Let $n$ be a positive integer. Then there exists a positive number $c_n$ such that for every $p\in(1,\be)$, $\a\in(0,1)$, $f\in\L_\a(\R^n)$, 
and for arbitrary $n$-tuples $(A_1,\cdots,A_n)$ and $(B_1,\cdots,B_n)$
of commuting self-adjoint operators with $A_j-B_j\in\bS_p$, $1\le j\le n$,
the operator $f(A_1,\cdots,A_n)-f(B_1,\cdots,B_n)$ belongs to $\bS_{p/\a}$ 
and the following inequality holds:
$$
\|f(A_1,\cdots,A_n)-f(B_1,\cdots,B_n)\|_{\bS_{p/\a}}
\le c_n\,(1-\alpha)^{-1}p^\a(p-1)^{-\a}\|f\|_{\L_\a}
\max_{1\le j\le n}\|A_j-B_j\|^\a_{\bS_p}.
$$
\end{thm}

Moreover, a stronger result holds: we can replace the $\bS_p$-norm with finite sums
$$
\left(\sum_{k=0}^m\big(s_k(T)\big)^p\right)^{1/p},
$$
where $T$ is a bounded linear operators on Hilbert space.

\begin{thm}
\label{Spnlki}
Let $n$ be a positive integer. Then there exists a positive number $c_n$ such that for every $m\in\Z_+$, $p\in(1,\be)$, $\a\in(0,1)$, and $f\in\L_\a(\R^n)$,
the following inequality holds:
\begin{align*}
\sum_{k=0}^m&\Big(s_k\big(f(A_1,\cdots,A_n)-f(B_1,\cdots,B_n)\big)\Big)^{p/\a}\\[.2cm]
&\le(c_n)^{p/\a}(1-\a)^{-p/\a}p^p(p-1)^{-p}\|f\|_{\L_\a}^{p/\a}
\max_{1\le j\le n}\sum_{k=0}^m\big(s_k(A_j-B_j)\big)^p,
\end{align*}
whenever $(A_1,\cdots,A_n)$ and $(B_1,\cdots,B_n)$ are $n$-tuples of commuting self-adjoint operators such the operators $A_j-B_j$ are bounded, $1\le j\le n$.
\end{thm}

Theorem \ref{Spnki} can be generalized to the much more general case of 
quasinormed ideals $\fI$ with upper Boyd index $\b_\fI$ less than 1. Recall that for such an ideal $\fI$, the number $\bs{C}_\fI$ is defined in Subsection 2.2.

\begin{thm}
\label{nideal}
Let $n$ be a positive integer. There exists a positive integer $c_n$ such that
for every $\a\in(0,1)$, $f\in\L_\a(\R^n)$, for an arbitrary quasinormed ideal $\fI$ with $\b_\fI<1$, and for arbitrary $n$-tuples $(A_1,\cdots,A_n)$ and $(B_1,\cdots,B_n)$
of commuting self-adjoint operators with $A_j-B_j\in\fI$, $1\le j\le n$,
the operator $\big|f(A_1,\cdots,A_n)-f(B_1,\cdots,B_n)\big|^{1/\a}$
belongs to $\fI$
and the following inequality holds:
\begin{align*}
\Big\|\big|f(A_1,\cdots,A_n)-f(B_1,\cdots,B_n)\big|^{1/\a}\Big\|_\fI\le
c_n^{1/\a}\bs{C}_\fI(1-\a)^{-1/\a}\|f\|_{\L_\a}^{1/\a}
\max_{1\le j\le n}\|A_j-B_j\|_{\fI}.
\end{align*}
\end{thm}

Note that Theorem \ref{Spnki} does not generalize to the case $p=1$. Indeed, it was established in Theorem 9.9 of \cite{AP3} that for $\a>0$, there exist self-adjoint operators $A$ and $B$ on Hilbert space such that $\rank(A-B)=1$, and a function $f\in\L_\a(\R)$ such that
$$
s_k\big(f(A)-f(B)\big)\ge(1+k)^{-\a},\quad k\ge0,
$$
and so $f(A)-f(B)\not\in\bS_{1/\a}$. 

However, for $p=1$ a slightly weaker result holds.

Recall that for $q>0$ the weak type space $\bS_{q,\be}$ consists of compact operators $T$ on Hilbert space, for which
$$
\|T\|_{\bS_{q,\be}}\df\sup_{k\ge0}s_k(T)(1+k)^{1/q}.
$$

\begin{thm}
\label{sltip}
Let $n$ be a positive integer. Then there is a positive number $c_n$
such that for every $\a\in(0,1)$, $f\in\L_\a(\R^n)$, and for arbitrary $n$-tuples $(A_1,\cdots,A_n)$ and $(B_1,\cdots,B_n)$
of commuting self-adjoint operators with $A_j-B_j\in\bS_1$, $1\le j\le n$, the 
operator $f(A_1,\cdots,A_n)-f(B_1,\cdots,B_n)$ belongs to $\bS_{\frac1\a,\be}$
and the following inequality holds:
$$
\big\|f(A_1,\cdots,A_n)-f(B_1,\cdots,B_n)\big\|_{\bS_{\frac1\a,\be}}\le 
c_n(1-\a)^{-1}\|f\|_{\L_\a}\max_{1\le j\le n}\|A_j-B_j\|^\a_{\bS_1}.
$$
\end{thm}

We conclude this section with the result that shows that if we replace in the statement of Theorem \ref{sltip} the condition $f\in\L_\a(\R^n)$ with the slightly stronger condition $f\in B_{\be,1}^\a$, then we can make the stronger conclusion that $f(A_1,\cdots,A_n)-f(B_1,\cdots,B_n)\in\bS_{1/\a}$.

\begin{thm}
\label{s1besa}
Let $n$ be a positive integer. Then there is a positive number $c_n$
such that for every $\a\in(0,1)$, for every function $f$ in the Besov class $B_{\be,1}^\a(\R^n)$, and for arbitrary $n$-tuples $(A_1,\cdots,A_n)$ and $(B_1,\cdots,B_n)$
of commuting self-adjoint operators with $A_j-B_j\in\bS_1$, $1\le j\le n$, the 
operator $f(A_1,\cdots,A_n)-f(B_1,\cdots,B_n)$ belongs to $\bS_{1/\a}$
and the following inequality holds:
$$
\big\|f(A_1,\cdots,A_n)-f(B_1,\cdots,B_n)\big\|_{\bS_{1/\a}}\le 
c_n(1-\a)^{-1}\|f\|_{B_{\be,1}^\a}\max_{1\le j\le n}\|A_j-B_j\|^\a_{\bS_1}.
$$
\end{thm}

\medskip

As we have already mentioned in the introduction to this section, the proofs of the results of this section are practically the same as the proofs of the corresponding results for $n=1$ (see \cite{AP3}) and $n=2$ (see \cite{APPS2}) once we use inequality \rf{siId} obtained in \S\,5.

\

\section{\bf Commutator and quasicommutator estimates}
\setcounter{equation}{0}
\label{ComQ}

\

In this final section of this paper we obtain norm (or ideal norms) estimates of {\it quasicommutators}
$f(A_1,\cdots,A_n)R-Rf(B_1,\cdots,B_n)$ in terms the quasicommutators \lb$A_jR-RB_j$, $1\le j\le n$, where $R$ is a bounded linear operator. In the special case $A_j=B_j$, $1\le j\le n$, we deal with commutators, while in the special case $R=I$ we have the problem of estimating 
perturbations $f(A_1,\cdots,A_n)-f(B_1,\cdots,B_n)$ that was considered in previous sections.

We need the following formula that expresses $f(A_1,\cdots,A_n)R-Rf(B_1,\cdots,B_n)$ in terms of double operator integrals.

\begin{thm}
\label{quasido}
Let $f$ be an $L^\be$ function on $\R^n$, let $R$ be a bounded linear operator, and let $(A_1,\cdots,A_n)$ and $(B_1,\cdots,B_n)$ be $n$-tuples of commuting self-adjoint operators such that the operators $A_jR-RB_j$ are bounded. 
Suppose that $\Psi_j$, $1\le j\le n$, are Schur multipliers in
 $\fM_{\R^n\!,\R^n}$ that satisfy formula {\em\rf{pred}}
Then 
$$
f(A_1,\cdots,A_n)R-Rf(B_1,\cdots,B_n)=\sum_{j=1}^n\,\,\,
\iint\limits_{\R^n\times\R^n}\Psi_j(x,y)\,dE_A(x)(A_jR-RB_j)\,dE_B(y).
$$
\end{thm}

\Pf The proof is exactly the same as the proof of Theorem \ref{OLcs}.
One has just replace $f(A_1,\cdots,A_n)-f(B_1,\cdots,B_n)$ with
$f(A_1,\cdots,A_n)R-Rf(B_1,\cdots,B_n)$ and replace $A_j-B_j$ with $A_jR-RB_j$. $\bl$

%As in the case of formula \rf{nintf}, we first assume that the operators $A_j$ and $B_j$ are bounded for $1\le j\le n$. 
%We have
%\begin{align*}
%\iint\limits_{\R^n\times\R^n}\Psi_j(x,y)\,dE_A(x)(A_jR&-RB_j)\,dE_B(y)=
%\iint\limits_{\R^n\times\R^n}\Psi_j(x,y)\,dE_A(x)A_jR\,dE_B(y)\\[.2cm]
%&-\iint\limits_{\R^n\times\R^n}\Psi_j(x,y)\,dE_A(x)RB_j\,dE_B(y)\\[.2cm]
%=&\iint\limits_{\R^n\times\R^n}\Psi_j(x,y)x_j\,dE_A(x)R\,dE_B(y)\\[.2cm]
%&-\iint\limits_{\R^n\times\R^n}\Psi_j(x,y)y_j\,dE_A(x)R\,dE_B(y)\\[.2cm]
%=&\iint\limits_{\R^n\times\R^n}\Psi_j(x,y)(x_j-y_j)\,dE_A(x)R\,dE_B(y).
%\end{align*}
%Thus it follows from \rf{pred} that
%\begin{align*}
%\sum_{j=1}^n\,\,\,
%\iint\limits_{\R^n\times\R^n}\Psi_j(x,y)&\,dE_A(x)(A_jR-RB_j)\,dE_B(y)\\[.2cm]
%=&
%\sum_{j=1}^n\,\,\,
%\iint\limits_{\R^n\times\R^n}\Psi_j(x,y)(x_j-y_j)\,dE_A(x)R\,dE_B(y)\\[.2cm]
%=&\iint\limits_{\R^n\times\R^n}\big(f(x)-f(y)\big)\,dE_A(x)R\,dE_B(y)\\[.2cm]
%\!\!\!\!\!=&\iint\limits_{\R^n\times\R^n}f(x)\,dE_A(x)R\,dE_B(y)-
%\iint\limits_{\R^n\times\R^n}f(y)\,dE_A(x)R\,dE_B(y)\\[.2cm]
%=&f(A_1,\cdots,A_n)R-Rf(B_1,\cdots,B_n).
%\end{align*}
%
%In the general case, when the operators $A_1,\cdots,A_n$ and $B_1,\cdots,B_n$ are not necessarily bounded we use the same approximation procedure as in the proof of Theorem 5.2 of \cite{APPS2}. $\bl$

Theorem \ref{quasido} allows us to obtain analogs of the results of the preceding sections for quasicommutators.  
In particular, functions of class $B_{\be,1}^1(\R^n)$ admit the estimate: 
$$
\Big\|f(A_1,\cdots,A_n)R-Rf(B_1,\cdots,B_n)\Big\|\le \const\|f\|_{B_{\be,1}^1}
\max_{1\le j\le n}\|A_jR-RB_j\|,
$$
whenever $(A_1,\cdots,A_n)$ and $(B_1,\cdots,B_n)$ and $R$ is a bounded linear operator
such that the operators $A_jR-RB_j$ are bounded.

The following result is an analog of Theorem \ref{amC}:

%\begin{thm}
%\label{comLi}
%Suppose that $R$ is a bounded linear operator, and $(A_1,\cdots,A_n)$ and $(B_1,\cdots,B_n)$ are $n$-tuples of commuting self-adjoint operators such that the operators $A_jR-RB_j$ are bounded. 
%
%If $f$ is an $L^\be$ function on $\R^n$ whose Fourier transform is supported in the ball \lb$\{\xi\in\R^n:~\|\xi\|\le\s\}$, then
%$$
%\Big\|f(A_1,\cdots,A_n)R-Rf(B_1,\cdots,B_n)\Big\|\le c_n\,\s\|f\|_{L^\be}
%\max_{1\le j\le n}\|A_jR-RB_j\|
%$$
%for some positive number $c_n$.
%
%If $f\in B_{\be,1}^1(\R^n)$, then $f$ 
%$$
%\Big\|f(A_1,\cdots,A_n)R-Rf(B_1,\cdots,B_n)\Big\|\le c_n\|f\|_{B_{\be,1}^1}
%\max_{1\le j\le n}\|A_jR-RB_j\|
%$$
%for some positive number $c_n$.
%\end{thm}

%Note that as in the case of perturbations (see Theorem \ref{yadn}), we can replace the operator norm in Theorem \ref{comLi} with the (quasi)norm of an arbitrary quasinormed ideal with majorization property.

If $f\in\L_\o(\R^n)$ and $(A_1,\cdots,A_n)$, $(B_1,\cdots,B_n)$ and $R$ are as above, then
\begin{align}
\label{o*}
\Big\|f(A_1,\cdots,A_n)R&-Rf(B_1,\cdots,B_n)\Big\|\nonumber\\[.2cm]
&\le\const\|f\|_{\L_\o}\|R\|
\max_{1\le j\le n}\o_*\left(\frac{\|A_jR-RB_j\|}{\|R\|}\right).
\end{align}

%\begin{thm}
%\label{Lomn}
%Let $n$ be a positive integer. There exists a positive number $c_n$ such that 
%for every modulus of continuity $\o$, for every $f\in\L_\o(\R^n)$, the following inequality holds:
%$$
%\Big\|f(A_1,\cdots,A_n)R-Rf(B_1,\cdots,B_n)\Big\|\le c_n\|f\|_{\L_\o}\|R\|
%\max_{1\le j\le n}\o_*\left(\frac{\|A_jR-RB_j\|}{\|R\|}\right),
%$$
%whenever $(A_1,\cdots,A_n)$ and $(B_1,\cdots,B_n)$ are $n$-tuples of commuting self-adjoint operators and $R$ is a bounded linear operator such that the operators $A_jR-RB_j$ are bounded, $1\le j\le n$.
%\end{thm}
%
%Recall that the function $\o_*$ is defined by \rf{o_*}.

In particular, in the case $\o(t)=t^\a$, $\a\in(0,1)$, the following holds
%\begin{thm}
%\label{LalCn}
%Let $n$ be a positive integer. There exists a positive number $c_n$ such that 
%for every $\a\in(0,1)$, for every $f\in\L_\a(\R^n)$, the following inequality holds:
\begin{align}
\label{-1}
\Big\|f(A_1,\cdots,A_n)R&-Rf(B_1,\cdots,B_n)\Big\|\nonumber\\[.2cm]
&\le\const(1-\a)^{-1}\|f\|_{\L_\a}
\Big(\max_{1\le j\le n}\|A_jR-RB_j\|\Big)^\a\|R\|^{1-\a}.
\end{align}
%whenever $(A_1,\cdots,A_n)$ and $(B_1,\cdots,B_n)$ are $n$-tuples of commuting self-adjoint operators and $R$ is a bounded linear operator such that the operators $A_jR-RB_j$ are bounded, $1\le j\le n$.
%\end{thm}

%We also state here an analog of Theorem \ref{Spnki} for quasicommutators.
%
%\begin{thm}
%\label{SpComnki}
%Let $n$ be a positive integer. There exists a positive number $c_n$ such that for every $p\in(1,\be)$, $\a\in(0,1)$, $f\in\L_\a(\R^n)$, 
%the following inequality holds:
%\begin{align*}
%\|f(A_1,\cdots,A_n)R&-Rf(B_1,\cdots,B_n)\|_{\bS_{p/\a}}\\[.2cm]
%&\le c_n\,(1-\alpha)^{-1}p^\a(p-1)^{-\a}\|f\|_{\L_\a}
%\max_{1\le j\le n}\|A_jR-RB_j\|^\a_{\bS_p}\|R\|^{1-\a},
%\end{align*}
%whenever $(A_1,\cdots,A_n)$ and $(B_1,\cdots,B_n)$ are $n$-tuples of commuting self-adjoint operators and $R$ is a bounded linear operator such that $A_jR-RB_j\in\bS_p$, $1\le j\le n$.
%\end{thm}

We can also state and prove analogs of Theorems \ref{Spnki}, \ref{Lippert}, \ref{Spnlki}, \ref{nideal}, \ref{sltip}, and \ref{s1besa} for quasicommutators. The proofs of these analogs are practically the same as the proofs of the corresponding results for functions of perturbed $n$-tuples. Namely, one has to replace on the left hand-side the operators$f(A_1,\cdots,A_n)-f(B_1,\cdots,B_n)$ with the operators $f(A_1,\cdots,A_n)R-Rf(B_1,\cdots,B_n)$ and on the right-hand side the operators $A_j-B_j$ with the operators $A_jR-RB_j$.

Note that in the special case $n=2$ the results listed above were obtained earlier in \cite{APPS2}. Indeed, in the case $n=2$ the problem of estimating $f(A_1,A_2)R-Rf(B_1,B_2)$ in terms of $A_1R-RB_1$ and $A_2R-RB_2$ for pairs of commuting self-adjoint operators $(A_1,A_2)$ and $(B_1,B_2)$ can be reformulated in terms of normal operators:

\medskip

{\it Let $N_1$ and $N_2$ be normal operators and let $R$ be a bounded operator such that the operators $N_1R-RN_2$ and $N_1^*R-RN_2^*$ are bounded. For a function $f$ on the complex plane, estimate the quasicommutator $f(N_1)R-Rf(N_2)$ in terms of $N_1R-RN_2$ and $N_1^*R-RN_2^*$}.

\medskip

However, in the case of normal operators it is also natural to consider another problem for quasicommutators:

\medskip

{\it Let $N_1$ and $N_2$ be normal operators and let $R$ be a bounded operator such that the operator $N_1R-RN_2$ is bounded. For a function $f$ on the complex plane, estimate the quasicommutator $f(N_1)R-Rf(N_2)$ in terms of $N_1R-RN_2$}.

\medskip

The latter problem was treated in detail in \cite{AP6}.

To complete the paper we compare the results obtained in \cite{APPS2} and \cite{AP6}. 
We assume that $N_1$ and $N_2$ are normal operators and $R$ is a bounded linear operator.
It follows from the results of \cite{APPS2} that if $f\in B_{\be,1}^1(\R^2)$, then
$$
\big\|f(N_1)R-Rf(N_2)\big\|\le\const\|f\|_{B_{\be,1}^1}
\max\big\{\|N_1R-RN_2\|,\|N_1^*R-RN_2^*\|\big\}.
$$
On the other hand, it was shown in \cite{JW} that if $f$ is a continuous function on $\C$ such that 
$$
\big\|f(N_1)R-Rf(N_2)\big\|\le\const\|N_1R-RN_2\|
$$
for all bounded normal operators $N_1$ and $N_2$ and all bounded linear operators $R$, then $f$ is a linear function, i.e., $f(z)=az+b$ for some $a,\,b\in\C$.

In the case of Schatten--von Neumann norms the situation is quite different. Indeed, it was shown in \cite{AD} and \cite{Sh} that for $1<p<\be$, the following inequality holds:
$$
\|N_1^*R-RN_2^*\|_{\bS_p}\le\const\|N_1R-RN_2\|_{\bS_p},
$$
and so when we estimate $\big\|f(N_1)R-Rf(N_2)\big\|_{\bS_{p/\a}}$ for $f\in\L_a(\R^2)$, we do not need $\|N_1^*R-RN_2^*\|_{\bS_p}$ on the right-hand-side.

As for estimating the operator norms of $f(N_1)R-Rf(N_2)$ for functions of class $\L_\a(\R^2)$, $0<\a<1$, surprisingly, the results of \cite{AP6} show that such norms can be estimated only in terms of $\|N_1R-RN_2$: 
$$
\big\|f(N_1)R-Rf(N_2)\big\|\le\const(1-\a)^{-2}\|f\|_{\L_\a}\|N_1R-RN_2\|^\a\|R\|^{1-\a}.
$$
Note however, that the results of \cite{APPS2} allows us to replace on the right hand-side the factor $(1-\a)^{-2}$ with $(1-\a)^{-1}$ if we replace $\|N_1R-RN_2\|$ with $\max\{\|N_1R-RN_2\|,\|N^*_1R-RN^*_2\|\}$ (compare with \rf{-1}).

In the case of an arbitrary modulus of continuity $\o$ it was shown in \cite{AP6} that for $f\in\L_\o(\R^2)$, 
$$
\big\|f(N_1)R-Rf(N_2)\big\|\le\const\,\o_{**}\left(\frac{\|N_1R-RN_2\|}{\|R\|}\right),
$$
where
$$
\o_{**}\df(\o_*)_*,\quad\mbox{i.e.,}\quad
\o_{**}(\d)=\d\int_\d^\be\frac{\o_*(\t)}{\t^2}\,d\t,\quad\d>0.
$$
Again, if we compare the above estimate with results of \cite{APPS2}, we see that if we replace on the right $\|N_1R-RN_2\|$ with $\max\{\|N_1R-RN_2\|,\|N^*_1R-RN^*_2\|\}$, we can replace $\o_{**}$ with $\o_*$ (compare with \rf{o*}).

\

\

\footnotesize
\noindent
\begin{tabular}{p{8.7cm}p{5cm}}
F.L. Nazarov &  V.V. Peller \\
Department of Mathematics & Department of Mathematics  \\
Kent State University  & Michigan State University\\
Kent, Ohio 44242  & East Lansing, Michigan 48824 \\
USA & USA
\end{tabular}

\end{document}